\documentclass[12pt]{amsart}
\usepackage{latexsym,enumerate}
\usepackage{amsmath,amsthm,amsopn,amstext,amscd,amsfonts,amssymb}
\usepackage{color}
\usepackage{amssymb}
\numberwithin{equation}{section}

\newtheorem{theorem}{Theorem}
\newtheorem{lemma}{Lemma}

\newtheorem{proposition}{Proposition}
\newtheorem{remark}{Remark}

\numberwithin{claim}{section}

\numberwithin{theorem}{section}
\numberwithin{corollary}{section}
\numberwithin{lemma}{section}
\numberwithin{definition}{section}
\numberwithin{proposition}{section}
\numberwithin{remark}{section}

\newcommand{\medint}{-\kern  -,375cm\int}
\newcommand{\dint}{\displaystyle\int}

\begin{document}
\title[]{On the reverse isoperimetric inequality in Gauss space}
\author{ F. Brock$^{1}$ - F. Chiacchio$^{2}$ }
\thanks{}
\date{}

\begin{abstract}
\noindent In this paper we investigate the reverse isoperimetric inequality
with respect to the Gaussian measure for convex sets in $\mathbb{R}^{2}$.
While the isoperimetric problem for the Gaussian measure is well understood, many relevant aspects of the reverse problem have not yet been investigated.
In particular, to the best of our knowledge, there seem to be no results on the shape that the isoperimetric set should take.
Here,  through of a local perturbation analysis, 
we show that smooth perimeter-maximizing sets have
locally flat boundaries. Additionally, we derive sharper perimeter bounds
than those previously known, particularly for specific classes of convex
sets such as the convex sets symmetric with respect to the axes. Finally,
for quadrilaterals with vertices on the coordinate axes, we prove that the
set maximizing the perimeter ``degenerates" into the x-axis, traversed twice.

\vspace{0.2 cm}

\noindent \textsl{Key words: Gauss space, reverse isoperimetric inequality,
convex sets }

\vspace{0.2 cm}

\noindent \textsl{Mathematics Subject Classification: 52A40, 49K15, 49Q10}
\end{abstract}

\maketitle

\setcounter{footnote}{1} \footnotetext{Martin-Luther-Universität Halle-Wittenberg, Landesstudienkolleg, 
06114 Halle, Paracelsus-Str. 22, 
Germany, e-mail: friedemann.brock@studienkolleg.uni-halle.de}

\setcounter{footnote}{2} \footnotetext{
Dipartimento di Matematica e Applicazioni ``R. Caccioppoli", Universit\`{a}
degli Studi di Napoli Federico II, Complesso Monte S. Angelo, via Cintia,
80126 Napoli, Italy, e-mail: francesco.chiacchio@unina.it}

\section{Introduction}

In this paper, we study a reverse isoperimetric inequality in the plane
with respect to the Gaussian measure. Clearly the problem is meaningful only
if some constraints are imposed on the class of sets considered. Here we
focus on convex sets in $\mathbb{R}^{2}$. This problem, which has attracted
the attention of many authors, see e.g. \cite{B,L, L2, N} and the references
therein, falls within the category of shape optimization problems
(for an overview of this research topic, see, e.g., the classic monographs 
\cite{H} and \cite{HP}).

The problem addressed here has significant applications. For instance, its importance in Computational Learning Theory is discussed in \cite{AOS}
and \cite{J}  (see also the references therein).

It is worth noticing that, while the isoperimetric problem with respect to
the Gaussian measure has been thoroughly studied and fully understood, see 
\cite{Bo, CK, CFMP, ST}, many aspects of the reverse problem remain unclear.

Before discussing how our results fit into the existing literature, let us
also mention that similar issues have been addressed for the unweighted
case, see, e.g., \cite{B2,CZP} and the references therein.

Let $\Omega $ be a convex set in $\mathbb{R}^{N}$, with $N\geq 2,$ one
defines its Gaussian measure and perimeter as follows 
\begin{equation*}
\gamma _{N}\left( \Omega \right) =\frac{1}{\left( 2\pi \right) ^{\frac{N}{2}
} }\int\limits_{\Omega }\exp \left( -\frac{\left\vert \mathbf{x}\right\vert
^{2}}{2}\right) d\mathbf{x}
\end{equation*}
and 
\begin{equation*}
P_{\gamma _{N}}\left( \Omega \right) =\frac{1}{\left( 2\pi \right) ^{\frac{N 
}{2}}}\int_{\partial \Omega }\exp \left( -\frac{\left\vert \mathbf{x}
\right\vert ^{2}}{2}\right) d\mathcal{H}^{N-1}.
\end{equation*}
In \cite{B} the author proves that 
\begin{equation}
c_{N}:=\sup \left\{ P_{\gamma _{N}}\left( \Omega \right) ,\text{ }\Omega 
\text{ convex subset of }\mathbb{R}^{N}\right\} \leq 4N^{\frac{1}{4}}.
\label{c_N}
\end{equation}
Nazarov, in \cite{N}, among other things, showed that the behavior of the
constant $c_{N}$ with respect to $N$, provided in the previous inequality,
is optimal. However, for any fixed dimension $N$, it seems that neither the
value of best constant nor the shape of the maximizing set is known.

We are now in position to state our main results. The first result, derived
through a local analysis and presented in Section 2, shows that a
perimeter-maximizing set within the considered class is locally flat. So it
provides an indication of the shape that the maximizing set should take.

\begin{theorem}
\label{Local}
Let $\Omega $ be a maximizer for the shape functional: 
\begin{equation*}
D\text{ convex subset of \ }\mathbb{R}^{2}\rightarrow P_{\gamma _{2}}\left(
D\right) .
\end{equation*}
Assume that $\exists P\in \partial \Omega :$ $B_{2\delta }(P) \cap \partial \Omega  $ $\in C^{2},$
for some $\delta >0.$ Then $B_{\delta }(P)\cap \partial \Omega $ is a
straight segment.
\end{theorem}

Note that the Theorem above does not address the existence of a maximizing
set. The following result addresses this issue this, but within a
significantly small class of convex sets $\widetilde{T}$, which is the class
of convex quadrilaterals with the four vertices lying in the four different
semiaxes. In Section 3, we prove the following

\begin{theorem}
\label{Quadrilaterals}
 The following sharp inequality holds 
\begin{equation}
\label{Th}
P_{\gamma _{2}}\left( T\right) 
<
 \sqrt{\frac{2}{\pi }}
\text{ \ for any }
T\in \widetilde{T}.  
\end{equation}
\end{theorem}
\noindent We will prove the optimality of the constant in (\ref{Th}) by showing that
\begin{equation*}
\lim_{n\rightarrow +\infty }P_{\gamma _{2}}\left( T_{n}\right)
= 
\sqrt{\frac{2}{\pi }},
\end{equation*}
where $\left\{ T_n \right\}_{n \in \mathbb{N}}$ is a sequence of rhombuses that ``converges'' to the $x$-axis.
Thus, losely speaking, the optimal set within the class $\tilde{T}$ is the x-axis traversed twice.

In the last section, we obtain a bound for $c_{2}$ sharper than the one
presented in (\ref{c_N}). On the other hand, our result holds true in $
\widetilde{C}$, which is a class of convex sets much larger than $\widetilde{T}$, but
still strictly contained within the class of all convex sets in $\mathbb{R}
^{2}$.

More precisely, the class $\widetilde{C}$ is made by all convex bounded sets
in $\mathbb{R}^{2}$ such that 
\begin{equation*}
C\in \widetilde{C}\Leftrightarrow C=\left\{ (x,y)\in \left[ a,b\right]
\times \mathbb{R}:-g(x)\leq y\leq f(x)\right\},
\end{equation*}
where $a<0<b$, $f$, $g:\left[ a,b\right] \rightarrow \left[ 0,+\infty \right) $
are two concave functions such that 
\begin{equation*}
f(a)=f(b)=g(a)=g(b)=0
\end{equation*}
and 
\begin{equation*}
\max_{\left[ a,b\right] }f(x)=f(0),\ \min_{\left[ a,b\right] }g(x)=g(0).
\end{equation*}

\noindent In the fourth and last section, we prove the following

\begin{theorem}
\label{C} For any $C\in \widetilde{C}$ it holds that 
\begin{equation}
P_{\gamma _{2}}\left( C\right) =\frac{1}{2\pi }\int_{\partial C}\exp \left( -
\frac{x^{2}+y^{2}}{2}\right) d\mathcal{H}^{1}\leq \frac{2}{\sqrt{\pi }}.
\label{TH}
\end{equation}
\end{theorem}

\bigskip

\section{A maximizer is locally flat}

\noindent 
The next two Lemmas will be used in the proof (obtained arguing by contradiction) of Theorem \ref{Local}.
\begin{lemma}
\label{Local1} Let $\Omega \subset \mathbb{R}^{2}$ be a maximizer, $
(x_{0},y_{0})\in \partial \Omega $ with $y_{0}>0$ such that for some $\delta
>0\,$it holds 
\begin{equation*}
\partial \Omega \cap \left\{ (x,y)\in \mathbb{R}^{2}:\left\vert
x-x_{0}\right\vert \leq \delta \right\} =\left\{ (x,f(x))\in \mathbb{R}
^{2}:\left\vert x-x_{0}\right\vert \leq \delta \right\} ,
\end{equation*}
where $f\in C^{2}\left( \left[ x_{0}-\delta ,x_{0}+\delta \right] ;\mathbb{R}
\right) ,$ $\ f^{\prime }(x_{0})=0,$ \ $f^{\prime \prime }(x_{0})<0.$ Then 
\begin{equation}
f^{\prime \prime }(x_{0})\leq -2y_{0}.  \label{f''<}
\end{equation}
\end{lemma}

\bigskip

\noindent \textbf{Proof of Lemma \ref{Local1}.} \ Assumption (\ref{f''<}) tells us that, if $t>0$ is small enough, say $t\leq t_0 $, then there are unique numbers  $x_{1}(t),x_{2}(t)$ such that
$$ 
x_1 (t) <x_0 <x_2 (t) , \ \ 
f(x_{1}(t))=f(x_{2}(t))=y_{0}-t \ \mbox{ and }\ 
\left\vert x_{i}(t)-x_{0}\right\vert <\delta ,\text{ }i=1,2.
$$
Further, for $i=1,2,$ we have 
\begin{equation}
\lim_{t\rightarrow 0}x_{i}(t)=x_{0}  \label{xi->0}
\end{equation}
and 
\begin{equation}
x_{i}^{\prime }(t)=-\frac{1}{f^{\prime }\left( x_{i}(t)\right) }.
\label{x'i}
\end{equation}
Let 
\begin{equation*}
\Omega \left( t\right) :=\Omega \cap \left\{ (x,y)\in \mathbb{R}
^{2}:y<y_{0}-t\right\} .
\end{equation*}
We have 
\begin{equation*}
2\pi \left[ P_{\gamma _{2}}(\Omega \left( t\right) )-P_{\gamma _{2}}(\Omega
) \right] =\int_{x_{1}(t)}^{x_{2}(t)}\left[ \exp \left( -\frac{
x^{2}+(y_{0}-t)^{2}}{2}\right) -\exp \left( -\frac{x^{2}+f^{\prime }(x)^{2}}{
2}\right) \sqrt{1+f^{\prime }(x)^{2}}\right] dx
\end{equation*}
and 
\begin{equation*}
2\pi \frac{d}{dt}P_{\gamma _{2}}(\Omega \left( t\right) )=-\frac{1}{
f^{\prime }\left( x_{2}(t)\right) }\exp \left( -\frac{
x_{2}(t)^{2}+(y_{0}-t)^{2}}{2}\right) \left( 1-\sqrt{1+f^{\prime }\left(
x_{2}(t)\right) ^{2}}\right)
\end{equation*}
\begin{equation*}
+\frac{1}{f^{\prime }\left( x_{1}(t)\right) }\exp \left( -\frac{
x_{1}(t)^{2}+(y_{0}-t)^{2}}{2}\right) \left( 1-\sqrt{1+f^{\prime }\left(
x_{1}(t)\right) ^{2}}\right)
\end{equation*}
\begin{equation}
+\int_{x_{1}(t)}^{x_{2}(t)}\exp \left( -\frac{x^{2}+(y_{0}-t)^{2}}{2}\right)
(y_{0}-t)dx.  \label{L41}
\end{equation}
Notice that 
\begin{equation*}
f^{\prime }(x_{1}(t))>0>f^{\prime }(x_{2}(t))
\end{equation*}
and, by (\ref{x'i}) and (\ref{xi->0}), 
\begin{equation*}
\lim_{t\rightarrow 0}\frac{f^{\prime }(x_{2}(t))}{f^{\prime }(x_{1}(t))}
=\lim_{t\rightarrow 0}\frac{f^{\prime \prime }(x_{2}(t))}{f^{\prime \prime
}(x_{1}(t))}\frac{x_{2}^{\prime }(t)}{x_{1}^{\prime }(t)}=\lim_{t\rightarrow
0}\frac{f^{\prime \prime }(x_{2}(t))}{f^{\prime \prime }(x_{1}(t))}\frac{
f^{\prime }\left( x_{1}(t)\right) }{f^{\prime }\left( x_{2}(t)\right) }
=\lim_{t\rightarrow 0}\frac{f^{\prime }\left( x_{1}(t)\right) }{f^{\prime
}\left( x_{2}(t)\right) }.
\end{equation*}
Therefore 
\begin{equation}
\lim_{t\rightarrow 0}\frac{f^{\prime }(x_{2}(t))}{f^{\prime }(x_{1}(t))}=-1.
\label{-1}
\end{equation}
From (\ref{L41}) and (\ref{-1}) we deduce that 
\begin{equation}
2\pi \lim_{t\rightarrow 0}\left[ \frac{1}{f^{\prime }\left( x_{1}(t)\right) }
\frac{d}{dt}P_{\gamma _{2}}(\Omega \left( t\right) )\right]
=\lim_{t\rightarrow 0}\frac{-1}{f^{\prime }\left( x_{2}(t)\right) }\frac{
\exp \left( -\frac{x_{2}(t)^{2}+(y_{0}-t)^{2}}{2}\right) }{f^{\prime }\left(
x_{1}(t)\right) }\left( 1-\sqrt{1+f^{\prime }\left( x_{2}(t)\right) ^{2}}
\right)  \label{DPer}
\end{equation}
\begin{equation*}
+\lim_{t\rightarrow 0}\frac{1}{f^{\prime }\left( x_{1}(t)\right) ^{2}}\exp
\left( -\frac{x_{1}(t)^{2}+(y_{0}-t)^{2}}{2}\right) \left( 1-\sqrt{
1+f^{\prime }\left( x_{1}(t)\right) ^{2}}\right)
\end{equation*}
\begin{equation*}
+\lim_{t\rightarrow 0}\frac{1}{f^{\prime }\left( x_{1}(t)\right) }
\int_{x_{1}(t)}^{x_{2}(t)}\exp \left( -\frac{x^{2}+(y_{0}-t)^{2}}{2}\right)
(y_{0}-t)dx=:L_{1}+L_{2}+L_{3.}
\end{equation*}
We have 
\begin{equation}
L_{1}=\lim_{t\rightarrow 0}\exp \left( -\frac{x_{2}(t)^{2}+(y_{0}-t)^{2}}{2}
\right) \lim_{t\rightarrow 0}\frac{-f^{\prime }\left( x_{2}(t)\right) }{
f^{\prime }\left( x_{1}(t)\right) }\lim_{t\rightarrow 0}\frac{1-\sqrt{
1+f^{\prime }\left( x_{2}(t)\right) ^{2}}}{f^{\prime }\left( x_{2}(t)\right)
^{2}}=-\frac{1}{2}\exp \left( -\frac{x_{0}^{2}+y_{0}^{2}}{2}\right) ,
\label{L1}
\end{equation}
\begin{equation}
L_{2}=\lim_{t\rightarrow 0}\exp \left( -\frac{x_{2}(t)^{2}+(y_{0}-t)^{2}}{2}
\right) \lim_{t\rightarrow 0}\frac{1-\sqrt{1+f^{\prime }\left(
x_{1}(t)\right) ^{2}}}{f^{\prime }\left( x_{1}(t)\right) ^{2}}=-\frac{1}{2}
\exp \left( -\frac{x_{0}^{2}+y_{0}^{2}}{2}\right)  \label{L2}
\end{equation}
and finally 
\begin{equation*}
L_{3}=\lim_{t\rightarrow 0}\frac{-f^{\prime }\left( x_{1}(t)\right) }{
f^{\prime \prime }\left( x_{1}(t)\right) }\left\{ \frac{-1}{f^{\prime
}\left( x_{2}(t)\right) }\exp \left( -\frac{x_{2}(t)^{2}+(y_{0}-t)^{2}}{2}
\right) (y_{0}-t)\right.
\end{equation*}
\begin{equation}
\left. +\frac{1}{f^{\prime }\left( x_{1}(t)\right) }\exp \left( -\frac{
x_{1}(t)^{2}+(y_{0}-t)^{2}}{2}\right) (y_{0}-t)\right\} =-\frac{2y_{0}}{
f^{\prime \prime }(x_{0})}\exp \left( -\frac{x_{0}^{2}+y_{0}^{2}}{2}\right)
\label{L3}
\end{equation}
From (\ref{DPer}), (\ref{L1}), (\ref{L2}) and (\ref{L3}) we deduce that 
\begin{equation*}
2\pi \lim_{t\rightarrow 0}\left[ \frac{1}{f^{\prime }\left( x_{1}(t)\right) }
\frac{d}{dt}P_{\gamma _{2}}(\Omega \left( t\right) )\right] =\left( -1-\frac{
2y_{0}}{f^{\prime \prime }(x_{0})}\right) \exp \left( -\frac{
x_{0}^{2}+y_{0}^{2}}{2}\right) .
\end{equation*}
Since $\Omega $ is a maximizer, the above limit must be nonpositive. We
conclude that 
\begin{equation*}
-1-\frac{2y_{0}}{f^{\prime \prime }(x_{0})}\leq 0
\end{equation*}
and (\ref{f''<}) immediately follows. \hfill $\square $

\bigskip

\begin{lemma}
\label{Local2}  
Let $\Omega $ as in Lemma \ref{Local1}, then 
\begin{equation}
f^{\prime \prime }(x_{0})=-y_{0}.  \label{f''=}
\end{equation}
\end{lemma}

\noindent \textbf{Proof of Lemma \ref{Local2}.}\ Let $h\in C^{2}\left( \left[
x_{0}-\delta ,x_{0}+\delta \right] ;\mathbb{R}\right) ,$ suppose that supp$
(h)\subset \left[ x_{0}-\delta ,x_{0}+\delta \right] ,$ and define 
\begin{equation*}
f_{t}(x):=f(x)+th(x),\text{ for }\left\vert x-x_{0}\right\vert \leq \delta .
\end{equation*}
Clearly for small enough $|t|$ - say $\left\vert t\right\vert \leq \delta $
- we have $f_{t}^{\prime \prime }(x)<0$. Finally, let 
\begin{equation*}
\Omega (t):=\Omega \cap \left\{ (x,y)\in \mathbb{R}^{2}:y<f_{t}(x),\text{ }
\left\vert x-x_{0}\right\vert <\delta \right\} .
\end{equation*}
We have 
\begin{eqnarray*}
0 &\geq & 2\pi \left[ P_{\gamma _{2}}(\Omega (t))-P_{\gamma _{2}}(\Omega )\right]
\\
&
=& \int_{x_{0}-\frac{\delta }{2}}^{x_{0}+\frac{\delta }{2}}\left( \exp \left(
- \frac{x^{2}+f_{t}(x)^{2}}{2}\right) \sqrt{1+f_{t}^{\prime }(x)^{2}}-\exp
\left( -\frac{x^{2}+f(x)^{2}}{2}\right) \sqrt{1+f^{\prime }(x)^{2}}\right) dx
\end{eqnarray*}
which implies 
\begin{equation*}
0=2\pi \left. \frac{d}{dt}P_{\gamma _{2}}(\Omega (t))\right\vert
_{t=0}=\int_{x_{0}-\frac{\delta }{2}}^{x_{0}+\frac{\delta }{2}}\exp \left( - 
\frac{x^{2}+f(x)^{2}}{2}\right) \left[ -fh\sqrt{1+f^{\prime }(x)^{2}}+\frac{
f^{\prime }h^{\prime }}{\sqrt{1+f^{\prime }(x)^{2}}}\right] dx.
\end{equation*}
Taking into account that $h\left( x_{0}\pm \frac{\delta }{2}\right) =0$
and integrating by parts we obtain 
\begin{equation*}
0= \int_{x_{0}-\frac{\delta }{2}}^{x_{0}+\frac{\delta }{2}}\exp \left( -\frac{
x^{2}+f(x)^{2}}{2}\right) \left[ -fh\sqrt{1+f^{\prime }(x)^{2}}+\left(
x+ff^{\prime }\right) \frac{f^{\prime }h}{\sqrt{1+f^{\prime }(x)^{2}}}
-\left( \frac{f^{\prime }}{\sqrt{1+f^{\prime }(x)^{2}}}\right) ^{\prime }h 
\right] dx.
\end{equation*}
Now the Fundamental Lemma of the Calculus of variations yields 
\begin{equation*}
-f\sqrt{1+f^{\prime }(x)^{2}}+\left( x+ff^{\prime }\right) \frac{f^{\prime } 
}{\sqrt{1+f^{\prime }(x)^{2}}}-\left( \frac{f^{\prime }}{\sqrt{1+f^{\prime
}(x)^{2}}}\right) ^{\prime }=0\quad \text{on }\left( x_{0}-\frac{\delta }{2}
,x_{0}+\frac{\delta }{2}\right)
\end{equation*}
which implies 
\begin{equation*}
-f\sqrt{1+f^{\prime }(x)^{2}}+\frac{xf^{\prime }+ff^{\prime 2}}{\sqrt{
1+f^{\prime }(x)^{2}}}-\frac{f^{\prime \prime }}{\sqrt{1+f^{\prime }(x)^{2}}}
+\frac{f^{\prime 2}f^{\prime \prime }}{\left( 1+f^{\prime }(x)^{2}\right) ^{ 
\frac{3}{2}}}=0
\end{equation*}
\begin{equation*}
\Rightarrow \frac{-f\left( 1+f^{\prime }{}^{2}\right) +xf^{\prime
}+ff^{\prime 2}}{\sqrt{1+f^{\prime }(x)^{2}}}-\frac{f^{\prime \prime }}{
\left( 1+f^{\prime }(x)^{2}\right) ^{\frac{3}{2}}}=0
\end{equation*}
\begin{equation*}
\Rightarrow \left( -f\left( 1+f^{\prime 2}\right) +xf^{\prime }+ff^{\prime
2}\right) \left( 1+f^{\prime }{}^{2}\right) -f^{\prime \prime }=0
\end{equation*}
\begin{equation*}
\Rightarrow \left( -f+xf^{\prime }\right) (1+f^{\prime 2})-f^{\prime \prime
}=0.
\end{equation*}
Since $y_{0}=f(x_{0})$ and $f^{\prime }(x_{0})=0$ we get 
\begin{equation*}
-y_{0}-f^{\prime \prime }(x_{0})=0,
\end{equation*}
which is equality (\ref{f''=}). \hfill $\square $

\bigskip

\noindent \textbf{Proof of Theorem \ref{Local}}
Assume that the
assertion is not true. Then there exists a point $P\in \partial \Omega $ and $
\delta >0$, such that $B_{2\delta }(P)\cap \partial \Omega \in C^{2}$ and 
 $B_{\delta }(P)\cap \Omega $ is strictly convex. Note also that the convexity of $\Omega $ ensures that $B_{2\delta }(P)\cap \Omega $
lies on one side of $B_{2\delta }(P)\cap \partial \Omega $ only. 
Let $g$ be
the tangent line to $\partial \Omega $ passing through $P$,
by the strict convexity of $B_{\delta }(P) \cap \Omega $
 we have that  $\overline{\Omega } \cap g = \{ P\} $. 
\\
Suppose for a moment that the origin $O$ and $\Omega $ are
on opposite sides of $g,$ or $O\in g.$ Let ${\bf N}$ be the normal to $g$ pointing away from $\Omega $. Shifting $\Omega $ a little bit in direction ${\bf N }$, we then get a convex set  with larger Gaussian perimeter than $\Omega $, contradicting the assumptions. Hence we must have that $O\notin g$ and $O$ and $\Omega $
lie on the same side of $g$. 
\\
This means that we can achieve by a rotation of the coordinate system that 
\begin{equation*}
P=(x_{0},y_{0}),\text{ \ } \mbox{ with }\ y_{0}>0,
\end{equation*}
\begin{equation*}
\overline{B_{\delta }(P)}\cap \partial \Omega =\left\{ (x,f(x)):\left\vert
x-x_{0}\right\vert <\delta \right\} ,
\end{equation*}
\begin{equation*}
f(x_{0})=y_{0},\text{ }f^{\prime }(x_{0})=0\text{ \ and }f^{\prime \prime
}(x)<0\text{ on }\left[ x_{0}-\delta ,x_{0}+\delta \right] .
\end{equation*}
Therefore 
\begin{equation*}
\Omega \cap \left\{ (x,y):y>f(x),\text{ }\left\vert x-x_{0}\right\vert
<\delta \right\} =\emptyset .
\end{equation*}
But then 
Lemma \ref{Local1} and \ref{Local2} would ensure that 
\begin{equation*}
f^{\prime \prime }(x_{0})\leq -2y_{0}
\end{equation*}
and 
\begin{equation*}
f^{\prime \prime }(x_{0})= -y_{0},
\end{equation*}
which is impossible. 
$\hfill \Box $ \\[0.1cm]

In view of Theorem \ref{Local}, it is natural to ask whether the Gaussian perimeter attains its maximum for regular polygons centered at the origin. However, the next result shows that the answer is negative.

\bigskip

\begin{proposition}
\label{n-gon}

 Let $T_{n}(r)$ be the regular 
$n$-gon ($n\geq 3$)
symmetric with respect to the origin, such that the maximum distance of its boundary from the origin is $r$. It holds that
\begin{equation*}
P_{\gamma }(T_{n}(r))
<
 \sqrt{\frac{2}{\pi }}
\text{ \ }\forall r>0
\,\,\,\,  \text{and} \,\,\,\,
\forall n\geq 3.
\end{equation*}
\end{proposition}

\noindent \textbf{Proof of Proposition \ref{n-gon}.}\  Let us first analyze the case $n\geq 5.$

\noindent It holds that 
\begin{equation*}
P_{\gamma }(T_{n}(r))=\frac{n}{\pi }\exp \left( -\frac{r^{2}\cos ^{2}\left( 
\frac{\pi }{n}\right) }{2}\right) \int_{0}^{r\sin \left( \frac{\pi }{n}
\right) }\exp \left( -\frac{\sigma ^{2}}{2}\right) d\sigma \leq 
 r\exp \left( - \frac{r^{2}\cos ^{2}\left( \frac{\pi }{n}\right) }{2}\right).
\end{equation*}
Therefore $\forall n\geq 5$ and $\forall r>0$ we have 
\begin{equation*}
P_{\gamma }(T_{n}(r)\leq \left. r\exp \left( -\frac{r^{2}\cos ^{2}\left( 
\frac{\pi }{n}\right) }{2}\right) \right\vert _{r=\frac{1}{\cos \left( \frac{
\pi }{n}\right) }}=\frac{1}{\cos \left( \frac{\pi }{n}\right) }\exp \left( -
\frac{1}{2}\right) \leq \frac{1}{\cos \left( \frac{\pi }{5}\right) }\exp
\left( -\frac{1}{2}\right) =0.7497...
\end{equation*}
We get the claim since
\begin{equation*}
\frac{1}{\cos \left( \frac{\pi }{5}\right) }\exp \left( -\frac{1}{2}\right) 
<
 \sqrt{\frac{2}{\pi }}
=0.7979...
\end{equation*}

\noindent The case $n=4$ is contained in Theorem \ref{Quadrilaterals}. 

\noindent It remains to treat the case $n=3.$ 

\noindent  We have
\begin{equation*}
P_{\gamma }(T_{3}(r))=\frac{3}{\pi }\exp \left( -\frac{r^{2}}{8}\right)
\int_{0}^{\frac{\sqrt{3}}{2}r}\exp \left( -\frac{\sigma ^{2}}{2}\right)
d\sigma 
\end{equation*}
and
\begin{equation*}
\lim_{r\rightarrow 0^{+}}P_{\gamma }(T_{3}(r))=\lim_{r\rightarrow +\infty
}P_{\gamma }(T_{3}(r))=0.
\end{equation*}
This implies that the function $P_{\gamma }(T_{3}(r))$ attains its maximum at some
 $ r \in \left( 0,+\infty \right) .$

 Computing the derivative, we obtain
\begin{equation*}
\frac{d}{dr}P_{\gamma }(T_{3}(r))
=
\frac{3}{2\pi }\exp \left( -\frac{r^{2}}{8}
\right) \left[ \sqrt{3}\exp \left( -\frac{3r^{2}}{8}\right) -\frac{r}{2}
\int_{0}^{\frac{\sqrt{3}}{2}r}\exp \left( -\frac{\sigma ^{2}}{2}\right)
d\sigma \right] .
\end{equation*}
Thus
\begin{equation*}
\frac{d}{dr}P_{\gamma }(T_{3}(r))\underset{(<)}{>}0\text{ \ }\Leftrightarrow 
\text{ \ }g(r):=\sqrt{3}\exp \left( -\frac{3r^{2}}{8}\right) -\frac{r}{2}
\int_{0}^{\frac{\sqrt{3}}{2}r}\exp \left( -\frac{\sigma ^{2}}{2}\right)
d\sigma \underset{(<)}{>}0.
\end{equation*}
Since 
\begin{equation*}
\lim_{r\rightarrow 0^{+}}g(r)=\sqrt{3}\text{ \ and \ }\lim_{r\rightarrow
+\infty }g(r)=-\infty ,
\end{equation*}
 the function $g(r)$ must must vanish at least once in 
$ \left( 0,+\infty \right) .$

Moreover, since $\sqrt{3}\exp \left( -\dfrac{r^{2}}{8}\right) $ is strictly
decreasing on $\left[ 0,+\infty \right) $ while $\dfrac{r}{2}\dint_{0}^{
\frac{\sqrt{3}}{2}r}\exp \left( -\frac{\sigma ^{2}}{2}\right) d\sigma $ is
strictly increasing on $\left[ 0,+\infty \right) $,
  it follows that there exists a unique
 $\widehat{r}\in \left( 0,+\infty \right) $ such that
\begin{equation*}
\left. \frac{d}{dr}P_{\gamma }(T_{3}(r))\right\vert _{r=\widehat{r}}=0
\,\,\,\,   \text{and} \,\,\,\,
P_{\gamma }(T_{3}(\widehat{r}))=\max_{\left[ 0,+\infty \right) }P_{\gamma
}(T_{3}(r)).
\end{equation*}
Since 
\begin{equation*}
g(1.49)>0\text{ \ and \ }g(1.50)<0,
\end{equation*}
we deduce that
\begin{equation*}
\widehat{r}\in \left( 1.49,1.50 \right).
\end{equation*}
A straightforward computation yields
\begin{equation*}
\left\vert P_{\gamma }(T_{3}(\widehat{r}))-P_{\gamma
}(T_{3}(1.49))\right\vert \leq \frac{1}{100}\max_{\left[ 1.49,1.50 \right]
}\left\vert \frac{d}{dr}P_{\gamma }(T_{3}(r))\right\vert <\frac{1}{100}.
\end{equation*}
We conclude by observing that
\begin{equation*}
\max_{\left[ 0,+\infty \right) }P_{\gamma }(T_{3}(r))=P_{\gamma }(T_{3}(
\widehat{r}))\leq P_{\gamma }(T_{3}(1.49))+\frac{1}{100}<0.7382
<
 \sqrt{\frac{2}{\pi }}
=0.7979...
\end{equation*}
This completes the proof of the Lemma.\hfill $\square $

\vspace{.2 cm}

\section{The case of convex quadrilaterals}

\vspace{.2 cm}

\noindent In this section we prove Theorem 1.2.

\vspace{.1 cm}

Let us start by
showing that inequality (\ref{Th}) is optimal. To this aim consider the sequence of
rhombi $\left\{ T_{n}\right\} _{n\in \mathbb{N}}$ with vertices in the
points $(n,0),$ $\left( 0,\dfrac{1}{n}\right) ,$ $(-n,0)$ and $\left( 0,- 
\dfrac{1}{n}\right) $ . One has 
\begin{eqnarray*}
\lim_{n\rightarrow +\infty }P_{\gamma _{2}}\left( T_{n}\right)
&=&\lim_{n\rightarrow +\infty }\frac{4}{2\pi }\int_{0}^{n}\exp \left( -\frac{
1}{2}\left( x^{2}+\left( -\frac{x}{n^{2}}+\frac{1}{n}\right) ^{2}\right)
\right) \sqrt{1+\frac{1}{n^{4}}}dx \\
&=&\frac{2}{\pi }\int_{0}^{\infty }\exp \left( -\frac{x^{2}}{2}\right) dx= 
\sqrt{\frac{2}{\pi }}.
\end{eqnarray*}
Hence the optimal set degenerates into the $x$-axis, traversed two times.

In order to prove our claim we will obtain the assertion by treating each side of $T$ in the same way. Or, equivalently, through suitable reflections
with respect to the coordinate axes, the quadrilateral $T$ can be replaced
by a rhombus with vertices on the axes and larger Gaussian perimeter. So it
is sufficient to estimate the Gaussian length of just one side of $T$. In
other words we have to prove that for any $a,L>0,$ setting 
\begin{equation*}
f(x)=a-\frac{a}{L}x,\text{ with }x\in \left( 0,L\right) ,
\end{equation*}
one has 
\begin{eqnarray*}
\frac{\sqrt{2\pi }}{2} &=&\int_{0}^{\infty }\exp \left( -\frac{t^{2}}{2}
\right) dt>\int_{0}^{L}\exp \left( -\frac{t^{2}+\left[ f(t)\right] ^{2}}{2}
\right) \sqrt{1+\left( \frac{d}{dt}f(t)\right) ^{2}}dt \\
&=&\sqrt{1+\frac{a^{2}}{L^{2}}}\int_{0}^{L}\exp \left( -\frac{1}{2}\left(
t^{2}+\left( a-\frac{a}{L}t\right) ^{2}\right) \right) dt.
\end{eqnarray*}
Equivalently we want to prove that for any $a,L>0$ it holds 
\begin{equation*}
I:=\sqrt{\frac{L^{2}+a^{2}}{L^{2}}}\exp \left( -\frac{a^{2}}{2}\right)
\int_{0}^{L}\exp \left( -\frac{1}{2}\left( \dfrac{L^{2}+a^{2}}{L^{2}}t^{2}-2
\dfrac{a^{2}}{L}t\right) \right) dt<\int_{0}^{\infty }\exp \left( -\frac{
t^{2}}{2}\right) dt.
\end{equation*}
Clearly we have 
\begin{equation*}
I=\sqrt{\frac{L^{2}+a^{2}}{L^{2}}}\exp \left( -\frac{a^{2}}{2}\right)
\int_{0}^{L}\exp \left( -\frac{1}{2}\left( \frac{\sqrt{L^{2}+a^{2}}}{L}t-
\dfrac{a^{2}}{\sqrt{L^{2}+a^{2}}}\right) ^{2}+\frac{1}{2}\dfrac{a^{4}}{
L^{2}+a^{2}}\right) dt.
\end{equation*}
Performing the change of variables 
\begin{equation*}
\sigma :=
\frac{\sqrt{L^{2}+a^{2}}}{L}
t-\dfrac{a^{2}}{\sqrt{L^{2}+a^{2}}} ,
\end{equation*}
we obtain 
\begin{equation*}
I=\exp 
\left[ 
-\frac{1}{2}
\left( a^{2}-\dfrac{a^{4}}{L^{2}+a^{2}}\right) 
\right] 
\int_{-\frac{a^{2}}{\sqrt{L^{2}+a^{2}}}}^{\frac{L^{2}}{\sqrt{
L^{2}+a^{2}}}}\exp \left( -\frac{\sigma ^{2}}{2}\right) d\sigma .
\end{equation*}
Now setting 
\begin{equation*}
\alpha :=\frac{a^{2}}{\sqrt{L^{2}+a^{2}}}>0\text{ \ and\ \ }\beta :=\frac{
L^{2}}{\sqrt{L^{2}+a^{2}}}>0
\end{equation*}
we obtain

\begin{equation}
I=I(\alpha ,\beta )=\exp \left( -\frac{\alpha \beta }{2}\right)
\int_{-\alpha }^{\beta }\exp \left( -\frac{\sigma ^{2}}{2}\right) d\sigma .
\label{I}
\end{equation}
Note that for all non-negative $\alpha $ and $\beta $ we have
\begin{eqnarray}
 & & I(\alpha ,\beta )=I(\beta ,\alpha ) \quad \mbox{and}  
 \label{simm}
\\
& &
\label{Ialpha0} 
I(\alpha ,0) = \int_0 ^{\alpha } \mbox{exp}\, \left(-\frac{\sigma ^2}{2} \right) \, d\sigma < \sqrt{\frac{\pi}{2}} .
\end{eqnarray}
\hspace*{0.5cm}With the new variables $\alpha $ and $\beta$, the claim of Theorem 1.2 can be written as
\begin{equation}
\sup_{\mathbb{R}_{++}^{2}}I(\alpha ,\beta )=\lim_{t\rightarrow +\infty
}I(t,0)=\int_{0}^{\infty }\exp \left( -\frac{t^{2}}{2}\right) dt=\sqrt{ \frac{
\pi }{2}},  
\label{claim}
\end{equation}
where 
\begin{equation*}
\mathbb{R}_{++}^{2}=\left\{ (\alpha ,\beta )\in \mathbb{R}^{2}:\alpha ,\beta
\geq 0\right\} .
\end{equation*}
The proof of (\ref{claim}) is rather lengthy, so we divide it into
several Lemmas.

\vspace{0.1 cm}

\begin{lemma}
\label{Lemma 3.1}
  If $\{(\alpha _n , \beta _n )\} \subset \mathbb{R} ^2 _{++} $ with $\alpha _n ^2 + \beta _n ^2 \to +\infty $, then
\begin{equation}
\label{lim1}
\limsup_{n\to \infty} I(\alpha_n , \beta _n ) \leq \sqrt{\frac{\pi}{2}} .
\end{equation}  
\end{lemma}
\noindent \textbf{Proof of Lemma \ref{Lemma 3.1}.}
Suppose that (\ref{lim1}) does not hold. Then there exist  a number $\delta >0$ and a sequence $\{ (\alpha _n , \beta _n ) \} \subset \mathbb{R} ^2 _{++} $ with $\alpha _n ^2 + \beta _n ^2 \to +\infty $, such that 
\begin{equation}
\label{ineq1}
I(\alpha _n , \beta _n )\geq \sqrt{\frac{\pi}{2}} +\delta  \quad \forall n\in \mathbb{N} .
\end{equation}
In view of (\ref{simm}) we may assume w.l.o.g. that $\beta _n \to +\infty $. Further, since
$$
I (\alpha _n , \beta _n ) \leq \mbox{exp}\, \left(-\frac{\alpha _n \beta _n }{2} \right) \cdot \sqrt{2\pi } ,
$$
(\ref{ineq1}) yields that the sequence $\{ \alpha _n \cdot \beta _n \} $ is bounded. In turn, this means that $\alpha _n \to 0$. Finally, we obtain
\begin{eqnarray*}
\sqrt{\frac{\pi}{2}} +\delta & \leq & \mbox{exp}\, \left( \frac{\alpha _n \beta _n }{2} \right) \cdot \left( \sqrt{ \frac{\pi}{2} } + \delta \right) 
\leq   \mbox{exp}\, \left( \frac{\alpha _n \beta _n }{2} \right) \cdot I(\alpha _n , \beta _n )
\\
 &=& \int_{-\alpha _n } ^{\beta _n } \mbox{exp}\, \left( -\frac{ \sigma ^2 }{2} \right) \, d\sigma \longrightarrow \sqrt{\frac{\pi}{2} } \ \mbox{ as }\ n\to \infty ,
\end{eqnarray*}  
a contradiction.
$\hfill \Box$

\vspace{0.2 cm}

Next we study some properties of the mapping 
$$
t\longmapsto te^{-t^2 /2} , \quad ( t\in (0,+\infty ).
$$
Since the mapping is strictly increasing for $0\leq t<1$ and strictly decreasing for $1<t<+\infty $ with $\lim_{t\to +\infty } te^{-t^2 /2} =0$, there is for every $x\in (0, 1] $ a unique number $y=:f(x) \in [1, +\infty )$ such that
$$
x\cdot e^{-x^2 /2 } = y\cdot e^{-y^2 /2} .
$$
Further, there holds $f(1)=1$ and $\lim_{x\searrow 0} f(x)= +\infty $.  

\bigskip

\begin{lemma}
\label{Lemma 3.2} 
 The following properties hold:
\begin{eqnarray}
\label{C1} 
& & f\in C^1 (0, 1] ,\\
\label{fprime1}
& & f' (x) = \frac{(1-x^2) \cdot f(x)}{(1-f^2 (x) )\cdot x} \ \mbox{ for } x\in (0,1) , 
\\ 
\label{fprime2} 
& &  f'(1)=-1 ,
\\ 
\label{xf} 
& & x\cdot f(x) <1 \ \mbox{ for }\ x\in (0,1),
\\
\label{limxf} 
& & \lim_{x\searrow 0} x\cdot f(x) =0,
\\
\label{x+f} 
& & x+ f(x) >2 \ \mbox{ for }\ x\in (0, 1 ).
\end{eqnarray}
\end{lemma}

\noindent \textbf{Proof of Lemma \ref{Lemma 3.2}.}
 (\ref{C1}) and (\ref{fprime1}) follow from differentiation of
\begin{equation}
\label{xrelationtof}
x\cdot e^{-x^2 /2} = f(x) \cdot e^{-f^2 (x)/2} .
\end{equation}
Then (\ref{fprime2}) follows from (\ref{fprime1}) by passing to the limit $x\nearrow 1$.
\\
Assume that (\ref{xf}) was not true. Since $f(1)=1$, there exists $x_0 \in (0,1)$ with $x_0 \cdot f(x_0 ) \geq 1 $ and such that 
$$
\left( x \cdot f(x)\right) ^{\prime} |_{x=x_0 } \leq 0.
$$
By (\ref{fprime2}) we obtain
\begin{eqnarray*}
0 & \geq  & f(x_0 ) + x_0 \cdot f' (x_0 ) = f(x_0 ) + \displaystyle{\frac{ (1-x_0 ^2 ) \cdot f(x_0 )}{1-f^2 (x_0 )} }
\\
& = & \frac{f (x_0 ) \cdot (x_0 ^2 + f^2 (x_0 ) -2 )}{f^2 (x_0 ) -1 } ,
\end{eqnarray*}
which implies 
$$
x_0 ^2 + f^2 (x_0 ) \geq 2 .
$$
But this is impossible, since
$$
2 \leq 2x_0 \cdot f(x_0 ) < x_0 ^2 + f^2 (x_0 ). 
$$
Furthermore we calculate, using (\ref{xf}),
\begin{eqnarray*}
 \lim_{x\searrow 0} x\cdot f(x) &=&\lim_{x\searrow 0} \frac{ f(x)}{1/x}  = \lim_{x\searrow 0} \frac{f' (x)}{-1/x^2 } 
 \\
 & = & \lim_{x\searrow 0} \frac{(1-x^2 ) \cdot x\cdot f(x) }{f^2 (x)-1} =0. 
\end{eqnarray*}
Finally, (\ref{xf}) and (\ref{fprime1}) yield
\begin{eqnarray*}
(x+ f(x))^{\prime } & = & 1+ f' (x) = 1+ \frac{ (1-x^2 ) \cdot f(x)}{ (1-f^2 (x) \cdot x} 
\\
&=& \frac{(x+f(x)) \cdot (1- xf(x))}{(1-f^2 (x)) \cdot x} <0.
\end{eqnarray*}
Now (\ref{x+f}) follows from this and from the fact that $1+ f(1)=2$. 
$\hfill \Box $

\bigskip

\noindent Next we study the function $I(\alpha , \beta )$ at critical points.

\bigskip

\begin{lemma}
\label{Lemma 3.3} 
 Let $(\alpha , \beta ) \in (0, +\infty ) ^2 $ such that
$\nabla I(\alpha , \beta )=0$.
Then 
\begin{equation}
\label{nablaI2} 
\int_{-\alpha } ^{\beta } \mbox{exp}\, \left( \frac{-\sigma ^2}{2} \right) \, d\sigma = \frac{2}{\alpha } \cdot \mbox{exp}\, \left( -\frac{\beta ^2}{2} \right) = \frac{2}{\beta } \cdot \mbox{exp}\, \left( -\frac{\alpha ^2}{2} \right) ,
\end{equation}
and if $0<\alpha \leq 1$, then 
\begin{eqnarray}
\label{falphabeta}
&& \beta =f(\alpha )\ \mbox{ and}
\\
\label{Ialphafalpha1}
&& I(\alpha, \beta )= I(\alpha , f(\alpha )) = 
\\
\nonumber
& & \frac{2}{f(\alpha )} \cdot \mbox{exp}\, \left( - \frac{ \alpha f(\alpha ) +\alpha ^2 }{2} \right) = \frac{2}{\sqrt{\alpha f(\alpha )} } \cdot \mbox{exp}\, \left( -\frac{(\alpha + f(\alpha ))^2}{4} \right)  .
\end{eqnarray}
\end{lemma}

\noindent \textbf{Proof of Lemma \ref{Lemma 3.3}.}
We have 
\begin{eqnarray*} 
0 &=& I_{\alpha } = \mbox{exp}\, \left( -\frac{\alpha \beta}{2} \right) \cdot \left( 
-\frac{\beta }{2} \cdot \int_{-\alpha } ^{\beta } \mbox{exp}\, \left( -\frac{\sigma ^2}{2} \right) \, d\sigma + \mbox{exp}\, \left( -\frac{\alpha ^2}{2} \right) \right)
,
\\
0 &=& I_{\beta } = \mbox{exp}\, \left( -\frac{\alpha \beta}{2} \right) \cdot \left( 
-\frac{\alpha }{2} \cdot \int_{-\alpha } ^{\beta } \mbox{exp}\, \left( -\frac{\sigma ^2}{2} \right) \, d\sigma + \mbox{exp}\, \left( -\frac{\beta ^2}{2} \right) \right)
,
\end{eqnarray*}
and (\ref{nablaI2})--(\ref{Ialphafalpha1}) follow. 
$\hfill \Box $

\bigskip

Define
\begin{equation}
\label{u1}
u(x) := \frac{2}{f(x)} \cdot \mbox{exp}\, \left( - \frac{xf(x)+ x^2 }{2} \right) , \quad (0<x\leq 1 ).
\end{equation}
Note that (\ref{xrelationtof}) also yields
\begin{equation}
\label{u2}
u(x) = \frac{2}{\sqrt{xf(x)} } \cdot \mbox{exp}\, \left( - \frac{(x+f(x))^2}{4} \right) , \quad (0<x\leq 1 ).
\end{equation}

\bigskip

\begin{lemma}
\label{Lemma 3.4} 
 It holds 
\begin{equation}
\label{ux} 
u(x) < \sqrt{\frac{\pi}{2}} \quad \forall x\in (0, 1] .
\end{equation}
\end{lemma}

\noindent \textbf{Proof of Lemma \ref{Lemma 3.4}.}
We have $\lim_{x\searrow 0} u(x) =0$ and
$$
0<u(1) = \frac{2}{e} <\sqrt{\frac{\pi}{2}} .
$$
Thus, in order to prove (\ref{ux}), it is sufficient to show that
\begin{equation}
\label{uprime} \mbox{if $u'(x_0 )=0$ for some $x_0 \in (0,1)$, then $u(x_0 ) <\displaystyle{\sqrt{\frac{\pi}{2}}}$.}
\end{equation}
Let $u'(x_0)=0$, ($x_0 \in (0,1)$). Then we evaluate :
$$
0= \frac{2}{f^2 (x_0 )} \cdot \mbox{exp}\,  \left( - \frac{x_0 f(x_0 )+ x_0 ^2}{2} \right) \cdot \left[ -f(x_0 ) \cdot \frac{f(x_0 ) + x_0 f'(x_0 ) +2x_0 }{2} -f' (x_0 ) \right] .
$$ 
This implies
$$
0= f(x_0 ) \cdot ( f(x_0 ) + x_0 f' (x_0 ) +2x_0 ) +2f' (x_0 ) ,
$$
and in view of (\ref{fprime1}) we obtain
\begin{eqnarray}
\nonumber
&& 0=
f^2 (x_0 ) +2x_0 f(x_0 ) + \frac{1-x^2 ) \cdot f(x_0 )}{(1-f^2 (x_0 )) \cdot x} \cdot (x_0 f(x_0 ) +2)  \ \Rightarrow
\\
\nonumber
 && 0= (f(x_0 ) +2x_0 ) ( x_0 -f^2 (x_0 ) x_0 ) + (x_0 f(x_0 ) +2) (1-x_0 ^2 ) \ \Rightarrow
\\
\nonumber
&& 2 (f(x_0 ) x_0 +1) = x_0 f(x_0 ) \cdot (x_0 + f(x_0 )) ^2 \ \Rightarrow
\\
\label{x+ftoxf}
&& \frac{(x_0 + f(x_0 ))^2 }{2} -1 = \frac{1}{x_0 f(x_0 )} . 
\end{eqnarray}
Now (\ref{u2}) and (\ref{x+ftoxf}) yield
\begin{equation}
\label{ux0}
u(x_0 ) = 2 \cdot \mbox{exp}\, \left( - \frac{(x_0 + f(x_0 ))^2 }{4} \right) \cdot \sqrt{ \frac{(x_0 + f(x_0))^2}{2} -1} \ .
\end{equation}
Finally we set
\begin{equation}
\label{vdef}
v(t):= 2e^{-t} \cdot \sqrt{2t-1}, \quad (t\geq 1).
\end{equation}
Since $v(1) =2/e$ and 
\begin{eqnarray*} 
v' (t) &=& -2e^{-t} \cdot \sqrt{2t-1}+ 2 e^{-t} \cdot \frac{1}{\sqrt{2t-1}} 
\\
&=& \frac{4e^{-t}}{\sqrt{2t-1}} \cdot (1-t) \leq 0 \quad \mbox{for }\ t\geq 1,
\end{eqnarray*}
it follows that
\begin{equation}
\label{v<}
v(t) \leq \frac{2}{e} \quad \mbox{for }\ t\geq 1.
\end{equation}
By (\ref{x+f}) we have 
$$
\frac{(x_0 + f(x_0 ))^2}{4} \geq 1.
$$  
Hence (\ref{v<}) yields
$$
u(x_0 ) = v\left( \frac{(x_0 + f(x_0 ))^2}{4} \right) \leq \frac{2}{e} < \sqrt{\frac{\pi}{2}} ,
$$
which is (\ref{uprime}). The Lemma is proved.
$\hfill \Box $

\bigskip

\begin{lemma}
\label{Lemma 3.5} 
 Let $(\alpha , \beta ) \in (0, +\infty )^2 $ such that $\nabla I(\alpha ,\beta )=0$. 
Then 
\begin{equation}
\label{crit<}
I(\alpha , \beta )< \sqrt{\frac{\pi}{2}} .
\end{equation}
\end{lemma}

\noindent \textbf{Proof of Lemma \ref{Lemma 3.5}.}
 By (\ref{Ialphafalpha1}) and (\ref{u1}) we have that
$$
I(\alpha , \beta ) = I(\alpha , f(\alpha )) = u(\alpha ),
$$
and the assertion follows from  Lemma \ref{Lemma 3.4}. 
$\hfill \Box $

\bigskip

\noindent \textbf{Proof of Theorem  \ref{Quadrilaterals}.}
The proof of (\ref{claim}), and, hence, of Theorem \ref
{Quadrilaterals}, is a direct consequence of 
(\ref{Ialpha0}), Lemma 3.1 and  Lemma 3.5.
$\hfill \Box $ 

\section{A reverse isoperimetric inequality in a suitable class of convex
sets}

\noindent \textbf{Proof of Theorem \ref{C}} Without loss of generality we
may assume that $f(x)$ and $g(x)$ are smooth and strictly concave functions.
Inequality (\ref{TH}) is an immediate consequence of the following upper
bound for $J,$ the Gaussian measure of the portion of $\partial C$ contained
in the first quadrant: 
\begin{equation}
J=\int_{0}^{b}\exp \left( -\frac{x^{2}+\left[ f(x)\right] ^{2}}{2}\right) 
\sqrt{1+\left( f^{\prime }(x)\right) ^{2}}dx\leq \sqrt{\pi }.  \label{M}
\end{equation}
Indeed by repeating the same considerations to the other portions of $
\partial C$ one gets 
\begin{equation*}
P_{\gamma _{2}}\left( C\right) \leq 4\frac{1}{2\pi }\sqrt{\pi }=\frac{2}{
\sqrt{\pi }},
\end{equation*}
which is (\ref{TH}). In the first quadrant we have that 
\begin{equation*}
y=f(x)\Leftrightarrow x=g(y),
\end{equation*}
\begin{equation*}
1+\left( f^{\prime }(x)\right) ^{2}\leq \left( 1-f^{\prime }(x)\right) ^{2}
\end{equation*}
and 
\begin{equation*}
\sqrt{1+\left( f^{\prime }(x)\right) ^{2}}\leq 1-f^{\prime }(x).
\end{equation*}
Hence 
\begin{equation*}
J\leq \int_{0}^{b}\exp \left( -\frac{x^{2}+\left[ f(x)\right] ^{2}}{2}
\right) dt+\int_{0}^{b}\exp \left( -\frac{x^{2}+\left[ f(x)\right] ^{2}}{2}
\right) \left( -f^{\prime }(x)\right) dt
\end{equation*}
\begin{equation*}
=\int_{0}^{b}\exp \left( -\frac{x^{2}+\left[ f(x)\right] ^{2}}{2}\right)
dt+\int_{0}^{f(0)}\exp \left( -\frac{\left[ g(y)\right] ^{2}+y^{2}}{2}
\right) dt.
\end{equation*}
Note that 
\begin{equation*}
x^{2}+\left( f(0)-\frac{f(0)}{b}x\right) ^{2}=x^{2}\left( 1+\left( \frac{f(0)
}{b}\right) ^{2}\right) -2\frac{\left[ f(0)\right] ^{2}}{b}x+\left[ f(0)
\right] ^{2}
\end{equation*}
\begin{equation*}
=\left( x\sqrt{1+\left( \frac{f(0)}{b}\right) ^{2}}-\frac{\dfrac{\left[ f(0)
\right] ^{2}}{b}}{\sqrt{1+\left( \dfrac{f(0)}{b}\right) ^{2}}}\right) ^{2}-
\frac{\dfrac{\left[ f(0)\right] ^{4}}{b^{2}}}{1+\left( \dfrac{f(0)}{b}
\right) ^{2}}+\left[ f(0)\right] ^{2}
\end{equation*}
\begin{equation*}
=\left( x\sqrt{1+\left( \frac{f(0)}{b}\right) ^{2}}-\frac{\left[ f(0)\right]
^{2}}{\sqrt{b^{2}+\left[ f(0)\right] ^{2}}}\right) ^{2}+\frac{b^{2}\left[
f(0)\right] ^{2}}{b^{2}+\left[ f(0)\right] ^{2}}.
\end{equation*}
Similarly 
\begin{equation*}
\left( b-\frac{b}{f(0)}y\right) ^{2}+y^{2}=\left( y\sqrt{1+\left( \frac{b}{
f(0)}\right) ^{2}}-\frac{b^{2}}{\sqrt{\left[ f(0)\right] ^{2}+b^{2}}}\right)
^{2}+\frac{b^{2}\left[ f(0)\right] ^{2}}{b^{2}+\left[ f(0)\right] ^{2}}.
\end{equation*}
Therefore we have 
\begin{equation}
J\leq \int_{0}^{b}\exp \left\{ -\frac{1}{2}\left[ \left( x\sqrt{1+\left( 
\frac{f(0)}{b}\right) ^{2}}-\frac{\left[ f(0)\right] ^{2}}{\sqrt{b^{2}+\left[
f(0)\right] ^{2}}}\right) ^{2}+\frac{b^{2}\left[ f(0)\right] ^{2}}{b^{2}+
\left[ f(0)\right] ^{2}}\right] \right\} dx  \label{I_1_1_Ineq}
\end{equation}
\begin{equation*}
+\int_{0}^{a}\exp \left\{ -\frac{1}{2}\left[ \left( y\sqrt{1+\left( \frac{b}{
f(0)}\right) ^{2}}-\frac{b^{2}}{\sqrt{\left[ f(0)\right] ^{2}+b^{2}}}\right)
^{2}+\frac{b^{2}\left[ f(0)\right] ^{2}}{b^{2}+\left[ f(0)\right] ^{2}}
\right] \right\} dy
\end{equation*}
Performing the change of variables 
\begin{equation*}
\xi =x\sqrt{1+\left( \frac{f(0)}{b}\right) ^{2}}-\frac{\left[ f(0)\right]
^{2}}{\sqrt{b^{2}+\left[ f(0)\right] ^{2}}}
\end{equation*}
in the first integral at the right hand side of (\ref{I_1_1_Ineq}) and 
\begin{equation*}
\eta =y\sqrt{1+\left( \frac{b}{f(0)}\right) ^{2}}-\frac{b^{2}}{\sqrt{\left[
f(0)\right] ^{2}+b^{2}}}
\end{equation*}
in the second we obtain 
\begin{equation*}
J\leq \frac{b}{\sqrt{b^{2}+\left( f(0)\right) ^{2}}}\exp \left( -\frac{1}{2}
\frac{b^{2}\left[ f(0)\right] ^{2}}{b^{2}+\left[ f(0)\right] ^{2}}\right) 
\displaystyle\int_{-\frac{\left[ f(0)\right] ^{2}}{\sqrt{b^{2}+\left[ f(0)
\right] ^{2}}}}^{\frac{b^{2}}{\sqrt{b^{2}+\left[ f(0)\right] ^{2}}}}\exp
\left( -\frac{\xi ^{2}}{2}\right) d\xi 
\end{equation*}
\begin{equation*}
+\frac{f(0)}{\sqrt{\left( f(0)\right) ^{2}+b^{2}}}\exp \left( -\frac{1}{2}
\frac{b^{2}\left[ f(0)\right] ^{2}}{b^{2}+\left[ f(0)\right] ^{2}}\right) 
\displaystyle\int_{-\frac{b^{2}}{\sqrt{\left[ f(0)\right] ^{2}+b^{2}}}}^{
\frac{\left[ f(0)\right] ^{2}}{\sqrt{\left[ f(0)\right] ^{2}+b^{2}}}}\exp
\left( -\frac{\eta ^{2}}{2}\right) d\eta 
\end{equation*}
\begin{equation}
=\frac{b+f(0)}{\sqrt{b^{2}+\left( f(0)\right) ^{2}}}\exp \left( -\frac{1}{2}
\frac{b^{2}\left[ f(0)\right] ^{2}}{b^{2}+\left[ f(0)\right] ^{2}}\right) 
\displaystyle\int_{-\frac{\left[ f(0)\right] ^{2}}{\sqrt{b^{2}+\left[ f(0)
\right] ^{2}}}}^{\frac{b^{2}}{\sqrt{b^{2}+\left[ f(0)\right] ^{2}}}}\exp
\left( -\frac{t^{2}}{2}\right) dt.  \label{I_1_2_Ineq}
\end{equation}
Now setting 
\begin{equation*}
\alpha :=\frac{\left[ f(0)\right] ^{2}}{\sqrt{b^{2}+\left[ f(0)\right] ^{2}}}
>0,\text{ \ }\beta :=\frac{b^{2}}{\sqrt{b^{2}+\left[ f(0)\right] ^{2}}}>0
\end{equation*}
we obtain 
\begin{equation*}
\alpha \left( \alpha +\beta \right) =\frac{a^{2}}{\sqrt{b^{2}+a^{2}}}\left( 
\frac{\left[ f(0)\right] ^{2}}{\sqrt{b^{2}+\left[ f(0)\right] ^{2}}}+\frac{
b^{2}}{\sqrt{b^{2}+\left[ f(0)\right] ^{2}}}\right) =\left[ f(0)\right] ^{2},
\end{equation*}
\begin{equation*}
\beta \left( \alpha +\beta \right) =\frac{b^{2}}{\sqrt{b^{2}+\left[ f(0)
\right] ^{2}}}\left( \frac{\left[ f(0)\right] ^{2}}{\sqrt{b^{2}+a^{2}}}+
\frac{L^{2}}{\sqrt{b^{2}+a^{2}}}\right) =b^{2},
\end{equation*}
and 
\begin{equation*}
\alpha +\beta =\frac{\left[ f(0)\right] ^{2}}{\sqrt{b^{2}+\left[ f(0)\right]
^{2}}}+\frac{b^{2}}{\sqrt{b^{2}+\left[ f(0)\right] ^{2}}}=\sqrt{b^{2}+\left[
f(0)\right] ^{2}}.
\end{equation*}
With this change of variables (\ref{I_1_2_Ineq}) becomes 
\begin{equation}
J\leq \frac{\sqrt{\alpha }+\sqrt{\beta }}{\sqrt{\alpha +\beta }}\exp \left( -
\frac{\alpha \beta }{2}\right) \int_{-\alpha }^{\beta }\exp \left( -\frac{
t^{2}}{2}\right) dt=\frac{\sqrt{\alpha }+\sqrt{\beta }}{\sqrt{\alpha +\beta }
}I\left( \alpha ,\beta \right) ,  \label{Phi}
\end{equation}
where $I\left( \alpha ,\beta \right) $ is the function defined in (\ref{I}).
Clearly one has
\begin{equation*}
\frac{\sqrt{\alpha }+\sqrt{\beta }}{\sqrt{\alpha +\beta }}\leq \sqrt{2}\text{
\ }\forall (\alpha ,\beta )\in \mathbb{R}_{++}^{2}.
\end{equation*}
This last inequality, together with (\ref{claim}), gives the claim. $\hfill
\Box $ 

\bigskip 

\begin{remark}
Let 
\begin{equation*}
\Phi (\alpha ,\beta ):=\frac{\sqrt{\alpha }+\sqrt{\beta }}{\sqrt{\alpha
+\beta }}I\left( \alpha ,\beta \right) .
\end{equation*}
Softwares, such as Mathematica, Matlab and Octave, assert that 
\begin{equation*}
\max_{(\alpha ,\beta )\in \mathbb{R}_{++}^{2}}\Phi (\alpha ,\beta )=\Phi
(0.8769,0.8769)=1.4950=:M.
\end{equation*}
This would provide the bound 
\begin{equation*}
P_{\gamma _{2}}\left( C\right) \leq 4\frac{1}{2\pi }M=0.9517, \,\, \forall C\in \widetilde{C},
\end{equation*}
which, in turn, improves a bit on the one contained in Theorem \ref{C}
\begin{equation*}
P_{\gamma _{2}}\left( C\right) \leq \frac{2}{\sqrt{\pi }}=1.1284,  \,\, \forall C\in \widetilde{C}.
\end{equation*}
\end{remark}

\bigskip


\noindent \textbf{Acknowledgements}
The research of Francesco Chiacchio was partially supported
 by 
\lq\lq Linear and Nonlinear PDE's: New directions and Applications\rq\rq 
  project, CUP E53D23018060001, - funded by European Union - Next Generation EU within the PRIN 2022 PNRR program and by Gruppo Nazionale per l'Analisi Matematica, la Probabilit\`a e le loro Applicazioni (GNAMPA) of Istituto Nazionale di Alta Matematica (INdAM). 
The second author also thanks the University of Leipzig for kind hospitality.

\bigskip

%
%

%

\end{document}